\documentclass[11pt]{article}
\usepackage{amsfonts, amsmath, amssymb, amsthm}

\PassOptionsToPackage{obeyspaces}{url}
\usepackage[colorlinks=true,citecolor=blue,urlcolor=blue,linkcolor=blue,bookmarksopen=true]{hyperref}
\usepackage{tikz}
\usepackage{algorithm}
\usepackage{algpseudocode}

\usepackage{etoolbox}
\patchcmd{\thebibliography}{\leftmargin\labelwidth}{\leftmargin\labelwidth\addtolength\itemsep{-0.1\baselineskip}}{}{}

\renewcommand{\cal}[1]{\mathcal{#1}}

\newcommand{\comp}{\operatorname{comp}}
\newcommand{\diam}{\operatorname{diam}}
\newcommand{\cent}{\operatorname{center}}
\newcommand{\treextend}{\textproc{TreeExtend}}
\newcommand{\compextend}{\textproc{CompExtend}}
\newcommand{\query}{\textbf{Query} }
\newcommand{\ram}{R}
\newcommand{\oram}{\widetilde{R}}

\author{Joseph Briggs\thanks{Department of Mathematical Sciences, Technion Israel Institute of Technology, Haifa, Israel. \texttt{briggs@campus.technion.ac.il}.} \and Christopher Cox\thanks{Department of Mathematical Sciences, Carnegie Mellon University, Pittsburgh, PA, USA. \texttt{cocox@andrew.cmu.edu}. Supported in part by U.S.\ taxpayers through NSF CAREER grant DMS-1555149.}}

\title{Restricted online Ramsey numbers of matchings and trees}
\date{\today}

\newtheorem{theorem}{Theorem}
\newtheorem{lemma}[theorem]{Lemma}
\newtheorem{corollary}[theorem]{Corollary}
\newtheorem{claim}[theorem]{Claim}

\newtheorem{defn}[theorem]{Definition}
\newtheorem{question}[theorem]{Question}
\newtheorem{remark}[theorem]{Remark}

\newenvironment{thmref}[1]
  {\renewcommand{\thetheorem}{\ref{#1}$'$}%
   \addtocounter{theorem}{-1}%
   \begin{theorem}}
  {\end{theorem}}
\newenvironment{thmreff}[1]
  {\renewcommand{\thetheorem}{\ref{#1}$''$}%
   \addtocounter{theorem}{-1}%
   \begin{theorem}}
  {\end{theorem}}

\begin{document}
\maketitle
\abstract{Consider a two-player game between players Builder and Painter. Painter begins the game by picking a coloring of the edges of $K_n$, which is hidden from Builder.
    In each round, Builder points to an edge and Painter reveals its color.
	Builder's goal is to locate a particular monochromatic structure in Painter's coloring by revealing the color of as few edges as possible.
    The fewest number of turns required for Builder to win this game is known as the restricted online Ramsey number.
	In this paper, we consider the situation where this ``particular monochromatic structure'' is a large matching or a large tree.
	We show that in any $t$-coloring of $E(K_n)$, Builder can locate a monochromatic matching on at least ${n-t+1\over t+1}$ edges by revealing at most $O(n\log t)$ edges.
	We show also that in any $3$-coloring of $E(K_n)$, Builder can locate a monochromatic tree on at least $n/2$ vertices by revealing at most $5n$ edges.
}

\section{Introduction}

For families of graphs $\cal G_1,\dots,\cal G_t$, the \emph{Ramsey number}, $\ram(\cal G_1,\dots,\cal G_t)$, is the least integer $n$ such that any $t$-coloring of $E(K_n)$ contains a copy of some $G_i\in\cal G_i$ in color $i$ for some $i\in[t]$.
If $\cal G_1=\cal G_2=\dots=\cal G_t=\cal G$, then we abbreviate $\ram_t(\cal G)=\ram(\cal G_1,\dots,\cal G_t)$.
Furthermore, if $\cal G_i=\{G_i\}$, then we write $\ram(G_1,\dots,G_t)=\ram(\cal G_1,\dots,\cal G_t)$.
Determining the growth of the Ramsey number $\ram_2(K_n)$ is a central, wide-open problem in extremal combinatorics.

Decades of study have spawned a myriad of generalizations and variants of the Ramsey numbers.
In this paper, we consider a variant called the \emph{restricted online Ramsey numbers}, which is defined through a two-player game between players Builder and Painter.

Let $\cal G_1,\dots,\cal G_t$ be families of graphs and fix $n\geq\ram(\cal G_1,\dots,\cal G_t)$.
Painter begins the game by picking a coloring of the edges of $K_n$, which is hidden from Builder.
In each round, Builder points to an edge of $K_n$ and Painter reveals the color of that edge.
Builder wins the game once she locates a copy of some $G_i\in\cal G_i$ in color $i$ of Painter's coloring.

Builder's questions are called \emph{queries} and edges that Builder has queried previously in the game are called \emph{exposed}.
The restricted online Ramsey number, $\oram(\cal G_1,\dots,\cal G_t;n)$, is the least integer $\ell$ for which Builder can guarantee winning the game by querying at most $\ell$ edges of $K_n$, regardless of Painter's strategy.

The restricted online Ramsey numbers are themselves a variant on the \emph{online Ramsey numbers}, which were introduced independently by Beck~\cite{B93} and Kurek and Ruci\'nski~\cite{KR05}.
The (unrestricted) online Ramsey numbers are defined through the same Builder--Painter game as above, but do not restrict the number of vertices in play; thus, $\oram(\cal G_1,\dots,\cal G_t;\infty)$ can be used to denote the online Ramsey number.
Restricted online Ramsey numbers were originally alluded to in~\cite{CFGH2018} and were recently studied in their own right by Gonzalez, He and Zheng~\cite{gonzalez2019upper}.

For any $n\geq N:=\ram(\cal G_1,\dots,\cal G_t)$, including $n=\infty$, we arrive at the following straightforward bounds:
\begin{equation}\label{eqn:easy}
    {N\over 2}\leq\oram(\cal G_1,\dots,\cal G_t;n)\leq{N\choose 2}.
\end{equation}
It is therefore natural to wonder whether $\oram(\cal G_1,\dots,\cal G_t;n)$ is linear in $N$, quadratic in $N$, or somewhere in between.

In the case of cliques, with $N=\ram_2(K_k)$, Conlon~\cite{conlononline} showed that $\oram_2(K_k;\infty)\leq(1-\epsilon)^k{N\choose 2}$ for some fixed $\epsilon>0$ and $k$ sufficiently large; however, in the restricted setting, the best known upper bound is $\oram_2(K_k;N)\leq{N\choose 2}-\Omega(N\log N)$, due to Gonzalez, He and Zheng~\cite{gonzalez2019upper}.
This suggests that, in general, determining restricted online Ramsey numbers is substantially more difficult than determining their unrestricted counterparts.
\medskip

In this paper, we focus on the restricted online Ramsey numbers of matchings and trees.
Let $rK_2$ denote the matching on $r$ edges and let $\cal T_n$ denote the family of all trees on $n$ vertices.
Recall the following well-known Ramsey-type results.

\begin{theorem}[Cockayne and Lorimer~\cite{CL75}]\label{thm:rammatching}
    For $t\geq 2$ and positive integers $r_1,\dots,r_t$,
    \[
        \ram(r_1K_2,\dots,r_tK_2)=\max_i r_i+1+\sum_{i=1}^t(r_i-1).
    \]
\end{theorem}

\begin{theorem}[Folklore]\label{thm:2ramtree}
    For any graph $G$, either $G$ or the complement of $G$ is connected.
    In other words, $\ram_2(\cal T_n)=n$.
\end{theorem}

\begin{theorem}[Gerencs\'er and Gy\'arf\'as~\cite{GG67}]\label{thm:ramtree}
    $\ram_3(\cal T_n)=2n-2$ if $n$ is even and $\ram_3(\cal T_n)=2n-1$ if $n$ is odd.
\end{theorem}

Observe that for any $r,t$, we trivially have $\oram_t(rK_2;\infty)=\oram_t(\cal T_r;\infty)=(t-1)(r-1)+1$; however, the restricted version of the problem is more challenging.
Despite this, we show that the restricted online Ramsey numbers for these graphs are indeed linear in their respective Ramsey numbers.

Here and throughout this paper, we use $\lg=\log_2$ and employ the convention that $\lg 0=0$ so that our results can be stated uniformly.

\begin{thmref}{thm:rammatching}\label{thm:matching}
	Fix $t\geq 2$ and let $r_1,r_2,\dots, r_t$ be positive integers.
    If $n\geq \ram(r_1K_2,\dots,r_tK_2)$, then
    \[
        \oram(r_1K_2,\dots,r_tK_2;n)\leq{2t-1+(t-3)\lg(t-2)\over t+1}n.
    \]
\end{thmref}
In particular, Builder can win the game using $n$ queries for two colors, ${5\over 4}n$ queries for three colors, and in general, Builder can win the game using $\bigl(2+\lg(t-2)\bigr)n$ queries for $t$ colors.

Turning now to trees, the following result was communicated to us by Micek and Pegden.
\begin{thmref}{thm:2ramtree}[Micek and Pegden~\cite{mp-personal}]\label{thm:2tree}
    If $n\geq 2$, then $\oram_2(\cal T_n;n)=2n-3$.
\end{thmref}

Since the above result is unpublished, we present a proof in Section~\ref{sec:tree2} for completeness.
We extend this result to $3$-colorings in Section~\ref{sec:tree}.
Define $k(n):={n\over 2}+1$ if $n\equiv 2\pmod 4$ and $k(n):=\bigl\lceil{n\over 2}\bigr\rceil$ otherwise.
Observe that Theorem~\ref{thm:ramtree} states that any $3$-coloring of $E(K_n)$ contains a monochromatic tree on $k(n)$ vertices.
\begin{thmref}{thm:ramtree}\label{thm:tree}
    If $n\geq 3$, then $\oram_3(\cal T_{k(n)};n)\leq 5(n-1)$.
\end{thmref}
The proofs of Theorems~\ref{thm:matching},~\ref{thm:2tree} and~\ref{thm:tree} do not rely on Theorems~\ref{thm:rammatching},~\ref{thm:2ramtree} and~\ref{thm:ramtree}, and thus provide self-contained proofs of the existence of the respective monochromatic graphs as well.
\medskip

With the exception of Theorem~\ref{thm:2tree}, our upper bounds are not necessarily tight.
Indeed, in Theorem~\ref{thm:matching} with $t=4$, it is actually possible for Builder to locate one of these monochromatic matchings by querying only ${7\over 5}n$ edges, as opposed to the claimed bound of ${8\over 5}n$ queries.
We do not know whether or not the bound in Theorem~\ref{thm:tree} is tight.

However, we can prove a different form of tightness by extending the restricted online Ramsey numbers to consider the Builder--Painter game on the edges of $K_n$ when $n<\ram(\cal G_1,\dots,\cal G_t)$.
We consider two natural extensions.
\medskip

\noindent\textbf{The locating game.} In this game, Builder must either locate one of the monochromatic graphs or determine that Painter's coloring \emph{cannot} contain any of these graphs.

Suppose that after $\ell$ queries, Builder has exposed color classes $C_1,\dots,C_t$.
Builder has won the game if either
\begin{enumerate}
    \item For some $i\in[t]$, $C_i$ contains a copy of some $G_i\in\cal G_i$, or
    \item\label{lg:xclude} For every $\chi\colon E(K_n)\to[t]$ with $C_i\subseteq\chi^{-1}(i)$, the coloring $\chi$ does \emph{not} contain a copy of any $G_i\in\cal G_i$ in color $i$ for any $i\in[t]$.
\end{enumerate}

Define $\oram(\cal G_1,\dots,\cal G_t;n)$ to be the smallest $\ell$ for which Builder can guarantee to win the locating game on $K_n$ by querying at most $\ell$ edges.
This is an immediate extension of the original definition of $\oram$ since Case~\ref{lg:xclude} can never occur if $n\geq\ram(\cal G_1,\dots,\cal G_t)$.
Observe, however, that the lower bound in \eqref{eqn:easy} does not necessarily hold when $n<\ram(\cal G_1,\dots,\cal G_t)$.
\medskip

\noindent\textbf{The cornering game.} In this game, Painter is required to guarantee the existence of one of the monochromatic graphs and Builder must simply determine which color contains it.

Call a coloring $\chi\colon E(K_n)\to[t]$ \emph{valid} if it contains a copy of some $G_i\in\cal G_i$ in color $i$ for some $i\in[t]$.
Suppose that after $\ell$ queries, Builder has exposed color classes $C_1,\dots,C_t$.
Builder has won the game if there is some fixed $c\in[t]$ such that every valid coloring $\chi$ with $C_i\subseteq\chi^{-1}(i)$ contains a copy of some $G_c\in\cal G_c$ in color $c$.
In other words, Builder has won the game if she has determined that Painter's coloring \emph{must} contain one of the monochromatic graphs in color $c$.

Define $\oram'(\cal G_1,\dots,\cal G_t;n)$ to be the smallest $\ell$ for which Builder can guarantee to win the cornering game on $K_n$ by querying at most $\ell$ edges.
It is not immediately clear that $\oram'$ is an extension of $\oram$, but we will show in Section~\ref{sec:remarks} that $\oram'(\cal G_1,\dots,\cal G_t;n)=\oram(\cal G_1,\dots,\cal G_t;n)$ whenever $n\geq\ram(\cal G_1,\dots,\cal G_t)$.

\medskip

Observe that we always have the inequality $\oram(\cal G_1,\dots,\cal G_t;n)\geq\oram'(\cal G_1,\dots,\cal G_t;n)$.
With these extensions in hand, we show the following.

\begin{thmreff}{thm:rammatching}\label{thm:matchingtight}
    Fix $t,r\geq 2$. If $n=\ram_t(rK_2)-1$, then $\oram_t(rK_2;n)={n\choose 2}$ and $\oram'_t(rK_2;n)\geq\Omega(n^2)$.
\end{thmreff}

\begin{thmreff}{thm:ramtree}\label{thm:treetight}
    If $n\geq 3$, then $\oram_3(\cal T_{k(n)+1};n)\geq\oram'_3(\cal T_{k(n)+1};n)\geq\Omega(n^2)$.
\end{thmreff}

This paper is organized as follows.
In Section~\ref{sec:matching}, we prove Theorems~\ref{thm:matching} and~\ref{thm:matchingtight}.
We then prove Theorem~\ref{thm:2tree} in Section~\ref{sec:tree2} and prove Theorems~\ref{thm:tree} and~\ref{thm:treetight} in Section~\ref{sec:tree}.
In each of these sections, we focus first on describing Builder's strategy and then describe Painter's strategy.
We conclude with a list of open questions in Section~\ref{sec:remarks}.

\section{Monochromatic matchings}\label{sec:matching}

\paragraph{Builder's strategy.}
We begin by presenting the key lemma which motivates Builder's strategy.
Throughout this paper, a \emph{forest} is assumed to have no isolated vertices, i.e.\ every connected component is a tree with at least one edge.
For a forest $F$, we denote the set of connected components of $F$ by $\comp(F)$.

\begin{defn}\label{defn:goodforest}
	Let $F$ be a forest and let $\chi\colon E(F)\to[t]$.
	$F$ is said to be a \emph{good forest} (with respect to $\chi$) if
	\begin{itemize}
		\item $\chi$ is a proper edge-coloring of $F$, and
		\item there is some color $c\in[t]$ such that every component of $F$ contains an edge of color $c$.
	\end{itemize}
\end{defn}

\begin{lemma}\label{lem:goodforest}
	Fix $t\geq 2$ and let $r_1,\dots,r_t$ be positive integers.
	Suppose that $F$ is a forest and $\chi\colon E(F)\to[t]$ is a $t$-coloring.
	If $F$ is a good forest with respect to $\chi$ and $|V(F)|\geq\max_i r_i+\sum_i (r_i-1)$, then, for some $i\in[t]$, $F$ contains a matching of size $r_i$ in color $i$.
\end{lemma}

\begin{proof}
	Denote by $m_i$ the number of edges of color $i$ in $F$, so $\sum_i m_i=e(F)$.
	Since $F$ is a good forest, $\chi$ is a proper coloring of $E(F)$, so the largest matching in color $i$ in $F$ has size precisely $m_i$.
	Hence, we need only show that $m_i\geq r_i$ for some $i\in[t]$.

	By assumption, there is some color $c\in[t]$ which appears in each component of $F$, so $m_c\geq|\comp(F)|$.
	Therefore,
	\[
		m_c+\sum_i m_i \geq |\comp(F)|+e(F)=|V(F)|\geq\max_i r_i+\sum_i (r_i-1),
	\]
	so the claim follows from the pigeonhole principle.
\end{proof}

With this in mind, Builder's strategy is to locate a good forest in Painter's coloring which covers all but at most one vertex.
The following lemma presents the main tool employed by Builder to accomplish this.

\begin{lemma}\label{lem:grow}
	Let $\chi\colon E(K_n)\to[t]$, and suppose that Builder has exposed all edges of some tree $T\subseteq E(K_n)$ on $m\geq 1$ edges.
	Suppose $T$ is properly edge-colored under $\chi$ and that $xy\in E(K_n)$ is an exposed edge completely disjoint from $T$ with $\chi(xy)\notin\chi(T):=\{\chi(f):f\in E(T)\}$.

	There exists a procedure $\treextend(\chi,T,xy)$ that, by querying at most $1+\bigl\lfloor\lg\bigl(\diam(T)-1\bigr)\bigr\rfloor$ extra edges, returns a tree $T^*$ with the following properties:
	\begin{enumerate}
		\item\label{item:vertex} $V(T^*)\subseteq V(T)\cup \{x,y\}$,
		\item\label{item:edge} $e(T^*)\in\{m+1,m+2\}$,
		\item\label{item:colors} $\chi(T^*)\supseteq \chi(T)\cup\{\chi(xy)\}$, and
		\item\label{item:propcol} $T^*$ is properly edge-colored under $\chi$.
	\end{enumerate}
\end{lemma}

\begin{proof} We first define $\treextend$ and then prove the claimed properties.

Recall that a vertex $v$ of a tree $T$ is called a \emph{center} if $v$ is at a distance at most $\bigl\lceil \diam(T)/2\bigr\rceil$ from every other vertex of $T$, where $\diam(T)$ is the diameter of $T$.
Note that there will be one such vertex when $\diam(T)$ is even, and two such vertices if $\diam(T)$ is odd.
We denote the center vertex of $T$ by $\cent(T)$, where an arbitrary choice is made if there are two such vertices.
Additionally, for an edge $xy\in E(T)$, define $T(x,y)$ to be the subtree of $T$ which is formed by rooting $T$ at $y$ and removing all descendants of $x$.
\pagebreak

\begin{algorithmic}[1]
\Procedure{TreeExtend}{$\chi,T,xy$}
	\State Fix any proper 2-coloring $\eta\colon V(T)\to \{x,y\}$ of $T$\Comment{For $v\in V(T)$, write $\eta_v=\eta(v)$}
	\State $T'\gets T$
	\State $v\gets\cent(T)$
	\Loop
		\State \query the edge $v\eta_v$
		\If{$\chi(v\eta_v)=\chi(xy)$}
			\State\Return $T+v\eta_v$ \label{algl:ret1}
		\ElsIf{$T+v\eta_v$ is properly colored}
			\State\Return $T+v\eta_v+xy$ \label{algl:ret2}
		\Else
			\State There is some edge $vv'\in E(T')$ with $\chi(vv')=\chi(v\eta_v)$
			\If{$v'$ is a leaf of $T$}
				\State\Return $T-vv'+v\eta_v+xy$ \label{algl:ret3}
			\ElsIf{$v'$ is a leaf of $T'$}
				\State\Return $T-vv'+v\eta_v+xy+v'\eta_{v'}$ \label{algl:ret4}
			\Else\label{algl:lastelse}
				\State $T'\gets T'(v,v')$
				\State $v\gets\cent(T'(v,v'))$
			\EndIf
		\EndIf
	\EndLoop
\EndProcedure
\end{algorithmic}

	We show first that $\treextend(\chi,T,xy)$ does in fact return $T^*$ and bound the number of queries made in the process.

	A new edge is queried only when reaching the beginning of the loop.
	If $\diam(T')\in\{1,2\}$, then the procedure will return $T^*$ before reaching Line~\ref{algl:lastelse}, thus requiring only one additional query.
	Furthermore, if we reach Line~\ref{algl:lastelse}, then we will have $\diam(T'(v,v'))\leq\bigl\lceil{\diam(T')\over 2}\bigr\rceil$ since $v=\cent(T')$.
	From this we conclude that $\treextend(\chi,T,e)$ returns $T^*$ by querying at most $1+\bigl\lfloor\lg\bigl(\diam(T)-1\bigr)\bigr\rfloor$ extra edges, recalling that $\lg 0=0$.

	We now verify the claimed properties of $T^*$.
	Set $c=\chi(xy)$ and consider the four situations in which $\treextend(\chi,T,xy)$ can return $T^*$.
	\begin{itemize}
		\item $T^*$ is returned on Line~\ref{algl:ret1}.
			Here $V(T^*)=V(T)\cup\{\eta_v\}\subseteq V(T)\cup \{x,y\}$ and $e(T^*)=m+1$.
			Additionally, $\chi(v\eta_v)=c$, so since $c\notin\chi(T)$, we know that $T^*$ is properly edge-colored and $\chi(T^*)=\chi(T)\cup\{c\}$.
			Finally, $T^*$ is in fact a tree since adding the edge $v\eta_v$ does not create a cycle.
		\item $T^*$ is returned on Line~\ref{algl:ret2}.
			Here $V(T^*)=V(T)\cup\{x,y\}$ and $e(T^*)=m+2$.
			Additionally, we know that $\chi(v\eta_v)\neq c$, so since $T+v\eta_v$ is properly edge-colored, we know that $T^*$ is also properly edge-colored and $\chi(T^*)=\chi(T)\cup\{\chi(v\eta_v),c\}\supseteq \chi(T)\cup\{c\}$.
			Finally, $T^*$ is a tree since we do not create a cycle upon adding the edges $v\eta_v$ and $xy$.
		\item $T^*$ is returned on Line~\ref{algl:ret3}.
			Here we have $V(T^*)=(V(T)\setminus\{v'\})\cup\{x,y\}$ and $e(T^*)=m+1$.
			Now, since we have $\chi(vv')=\chi(v\eta_v)$ for some $vv'\in E(T')$, we know that, since $T'$ is a subtree of $T$, $\chi(vv')=\chi(v\eta_v)\neq c$ since $c\notin \chi(T)$.
			Thus, $T^*$ is properly edge-colored since we removed the edge $vv'$, and also $\chi(T^*)=\chi(T)\cup\{c\}$ since $\chi(v\eta_v)=\chi(vv')$.
			Finally, $T^*$ is a tree since $v'$ is a leaf of $T$, so $T-vv'$ is still at tree, and then adding the edges $v\eta_v$ and $xy$ does not create a cycle.
		\item $T^*$ is returned on Line~\ref{algl:ret4}.
			Here we have $V(T^*)=V(T)\cup\{x,y\}$ and $e(T^*)=m+2$.
			Now, since $v'$ is a leaf of $T'$ but not a leaf of $T$, this means that we must have previously queried the edge $v'\eta_{v'}$ and found that $\chi(v'\eta_{v'})=\chi(vv')$.
			Since $\eta$ is a proper 2-coloring of $V(T)$ and $vv'\in E(T)$, we must have $\eta_v\neq\eta_{v'}$, so $T^*$ is indeed properly edge-colored, and $\chi(T^*)=\chi(T)\cup\{c\}$.
			Finally, $T^*$ is again a tree since we added the path $v\eta_v\eta_{v'}v'$ and removed the edge $vv'$.\qedhere
	\end{itemize}
\end{proof}

We state a consequence for later reference.
\begin{corollary}\label{cor:allowedvalues}
	Let $\chi,T,xy$ be as in Lemma~\ref{lem:grow} and set $T^*=\treextend(\chi,T,xy)$.
	\begin{itemize}
		\item If $e(T)=1$, then $T^*$ either has $|\chi(T^*)|=2$ and is a path with $2$ edges, or has $|\chi(T^*)|=3$ and is a path with $3$ edges.
		\item If $e(T)=2$, then $T^*$ either has $|\chi(T^*)|\geq 4$, or has $|\chi(T^*)|=3$ and is a star with $3$ edges.
	\end{itemize}
\end{corollary}

We now have all of the necessary tools to describe Builder's full strategy and prove Theorem~\ref{thm:matching}.

\begin{proof}[Proof of Theorem~\ref{thm:matching}]
Builder maintains and grows a good forest $F$.
While there are still at least 2 vertices $x,y$ uncovered by $F$, Builder queries the edge $xy$.
If its color $\chi(xy)$ is already present among all components of $F$, Builder adds this edge (as a 2-vertex component) to $F$, and repeats.
Otherwise, there is some connected component $T\in\comp(F)$ which does not have an edge of color $\chi(xy)$.

Here, Builder uses $\treextend(\chi,T,xy)$ to return a tree $T^*$ and replaces $F$ by $F-T+T^*$.
By Lemma~\ref{lem:grow}, $F-T+T^*$ is also a good forest and covers at least one more vertex than $F$, so this process must eventually terminate.
Furthermore, the process terminates when $F$ covers all but at most one vertex of $K_n$, and thus, by Lemma~\ref{lem:goodforest}, Builder has located a copy of $r_iK_2$ in color $i$ for some $i\in[t]$.
Let $F^*$ be the forest found by Builder.

We now count the total number of queries used to build $F^*$.
For integers $m,k\geq 2$, define $q(m,k):=2k-1+(k-3)\lg(m-2)$ and define also $q(1,1):=1$.
For $T\in\comp(F^*)$, let $q(T)$ denote the number of queries used by Builder to construct $T$.
\begin{claim}\label{claim:boundingtime}
	If $e(T)=m$ and $|\chi(T)|=k$, then $q(T)\leq q(m,k)$.
\end{claim}
\begin{proof}
	If $k=1$, then also $e(T)=1$, so certainly $q(T)=1=q(1,1)$.
	Thus suppose $k\geq 2$.

	Since $\treextend$ always appends at least one new color to a tree, we see that, for some $\ell\leq k$, there were trees $T_1,\dots,T_\ell$ and edges $e_1,\dots,e_\ell$ with $T_1$ being a single edge, $T_\ell=T$ and $T_{i+1}=\treextend(\chi,T_i,e_i)$ for all $i\in[\ell-1]$.
	Certainly $q(T_1)=1$, and, by Lemma~\ref{lem:grow}, if $d_i:=\diam(T_i)$, then
	\[
		q(T_{i+1})\leq q(T_i)+1+\bigl(1+\lfloor\lg(d_i-1)\rfloor\bigr),
	\]
	where the extra $+1$ comes from querying the edge $e_i$.
	Therefore,
	\[
		q(T)\leq 1+\sum_{i=1}^{\ell-1}\bigl(2+\lfloor\lg(d_i-1)\rfloor\bigr)=2\ell-1+\sum_{i=2}^{\ell-1}\lfloor\lg(d_i-1)\rfloor.
	\]
	By Lemma~\ref{lem:grow}, we know that $e(T_{i+1})\in\{e(T_i)+1,e(T_i)+2\}$, so certainly $d_i\leq e(T_i)\leq m-1$ for all $i\in[\ell-1]$.

	We now break into two cases based on $T_2$:
	\begin{itemize}
		\item $|\chi(T_2)|=2$: Since $T_2=\treextend(\chi,T_1,e_1)$, we know that $T_2$ must be a path on $2$ edges by Corollary~\ref{cor:allowedvalues}; thus $d_2=2$.
			As such,
			\[
				q(T)\leq 2\ell-1+\sum_{i=3}^{\ell-1}\lfloor\lg(d_i-1)\rfloor\leq 2k-1+(k-3)\lg(m-2)=q(m,k).
			\]
		\item $|\chi(T_2)|\geq 3$: Here, again by Corollary~\ref{cor:allowedvalues}, we must actually have $|\chi(T_2)|=3$ and $T_2$ is a path on $3$ edges, so $d_2=3$.
			Additionally, in this situation, we must also have $\ell\leq k-1$, and so we bound
			\[
				q(T)\leq 2\ell+\sum_{i=3}^{\ell-1}\lfloor\lg(d_i-1)\rfloor\leq 2k-2+(k-4)\lg(m-2)\leq q(m,k).\qedhere
			\]
	\end{itemize}
\end{proof}

For positive integers $m,k$, let $F_{m,k}$ denote the forest formed by all trees $T\in\comp(F^*)$ with $e(T)=m$ and $|\chi(T)|=k$.
By Lemma~\ref{lem:grow} and Corollary~\ref{cor:allowedvalues}, the only values of $(m,k)$ for which $F_{m,k}$ can be nonempty are: $(1,1)$, $(2,2)$, $(3,3)$ and $(m,k)$ where $4\leq k\leq m\leq 2k-3$.

For $4\leq k\leq m\leq 2k-3$, a quick calculation shows that
\[
	{q(m,k)\over m+1}\geq{q(m+1,k)\over m+2}\implies \max_{m:k\leq m\leq 2k-3}{q(m,k)\over m+1}={q(k,k)\over k+1}.
\]
Finally, another short calculation yields ${q(k,k)\over k+1}\leq{q(t,t)\over t+1}$ for all $k\in[t]$.

Thus, by Claim~\ref{claim:boundingtime}, we find that the total number of queries used to locate $F^*$, and thus the monochromatic matching, is bounded above by
\begin{align*}
	\sum_{T\in\comp(F^*)}q(T) &\leq  \sum_{m,k}|\comp(F_{m,k})|q(m,k)=\sum_{m,k}|V(F_{m,k})|{q(m,k)\over m+1}\\
				  &\leq {q(t,t)\over t+1}\sum_{m,k}|V(F_{m,k})|\leq {1\over t+1}\bigl(2t-1+(t-3)\lg(t-2)\bigr)n.\qedhere
\end{align*}
\end{proof}

\begin{remark}\label{rm:nottight}
	For any fixed $t\geq 4$, one can embark on a more sensitive analysis to improve the upper bound of ${1\over t+1}\bigl(2t-1+(t-3)\lg(t-2)\bigr)n$.
	For example, when $t=4$, we can arrive at an upper bound of ${7\over 5}n$ as opposed to ${8\over 5}n$ by working through Claim~\ref{claim:boundingtime} more carefully.
	However, as $t$ grows, it becomes increasingly difficult to carry out such an analysis.
\end{remark}
\begin{remark}\label{rm:constimprove}
	We can improve slightly the upper bound of $n$ queries when $t=2$.
	Indeed, with the strategy described above, the number of queries is at most ${1\over 2}|V(F_{1,1})|+|V(F_{2,2})|$, which is equal to $n$ only when $|V(F_{2,2})|=n$.
	However, $|V(F_{2,2})|$ is a multiple of $3$, so the bound can be improved to $n-1$ whenever $3\nmid n$.

	For $t=3$, similar reasoning will improve also the upper bound of ${5\over 4}n$ by an additive constant depending on the value of $n\pmod 4$.
\end{remark}

\paragraph{Painter's strategy.}
We now describe Painter's strategy in both the locating and cornering games for hiding a monochromatic $rK_2$ from Builder when $n=\ram_t(rK_2)-1=(t+1)r-t$.

\begin{proof}[Proof of Theorem~\ref{thm:matchingtight}]
	Consider partitioning $V(K_n)=V_1\sqcup\dots\sqcup V_t$ where $|V_1|=2r-1$ and $|V_i|=r-1$ for all $i\geq 2$.
	Let $\chi\colon E(K_n)\to[t]$ be given by $\chi(xy)=\max\{i:V_i\cap\{x,y\}\neq\varnothing\}$.

	Certainly $\chi$ does not contain a monochromatic $rK_2$.
	For an edge $e\in E(K_n)$ and a color $c\in[t]$, let $\chi_{e,c}$ denote the coloring of $E(K_n)$ obtained by giving $e$ color $c$ and coloring the rest of the edges as in $\chi$.
	Notice that if $e$ is not completely contained in $V_1$, then $\chi_{e,1}$ has a monochromatic $rK_2$ in color $1$, and if $e$ is completely disjoint from $V_c$ for some $c\geq 2$, then $\chi_{e,c}$ has a monochromatic $rK_2$ in color $c$.

    Consider first the locating game wherein Painter is not required to guarantee the existence of a monochromatic $rK_2$.
	As Builder queries edges, Painter colors the edge as in $\chi$ until there is only one unexposed edge, call it $e$.
	Certainly there is some $c\in[t]$ such that $\chi_{e,c}$ has an $rK_2$ in color $c$, so since $\chi$ does not have any monochromatic $rK_2$, Builder must query all ${n\choose 2}$ edges of $K_n$ in order to determine whether or not Painter's coloring has a monochromatic matching of size $r$.

    Consider now the cornering game wherein Painter is required to guarantee the existence of a monochromatic $rK_2$ and Builder needs only determine which color has said matching.
	Again, Painter will color the edges that Builder queries as in $\chi$ until the very last edge, which she then gives a color which will form a monochromatic $rK_2$.
	However, against this strategy, Builder can sometimes deduce which color will have this matching before reaching the very last edge.

	If $t=2$, then Painter's coloring has an $rK_2$ in color $2$ if and only if some edge in $V_1$ gets color $2$ (and otherwise the coloring must have an $rK_2$ in color $1$).
	Therefore, to determine which color contains the matching, Builder must either query every edge in $V_1$ or query every edge not completely contained in $V_1$,
	Hence, Builder must query at least
\[
	\min\biggl\{{|V_1|\choose 2},{n\choose 2}-{|V_1|\choose 2}\biggr\}={1\over 9}(2n+1)(n-1)
\]
	edges.

    For $t=3$, Builder must query all edges \emph{not} in $E[V_2,V_3]$ (the edges with one vertex in $V_2$ and the other in $V_3$); otherwise, there is some unexposed edge which is either completely contained in $V_1$ or meets $V_1$ in only one vertex and is disjoint from either $V_2$ or $V_3$.
	Therefore, Builder must query at least
	\[
		{n\choose 2}-|V_2||V_3|={1\over 16}(7n+1)(n-1)
	\]
	edges.

	Lastly, for $t\geq4$, since every edge is disjoint from at least two of $V_1,V_2,V_3,V_4$, Builder must query all ${n\choose 2}$ edges to determine which color contains the matching.
\end{proof}

\section{Monochromatic trees}
\subsection{Two colors}\label{sec:tree2}

Here we give the proof of Theorem~\ref{thm:2tree}, which was communicated to us by Micek and Pegden~\cite{mp-personal}.

\paragraph{Builder's strategy.}
Let $\chi\colon E(K_n)\to[2]$.
Builder begins by choosing $v\in V(K_n)$ arbitrarily and queries all $n-1$ edges incident to $v$.
For $i\in\{1,2\}$, set $C_i^0=\{v\}\cup\{u\in V(K_n):\chi(uv)=i\}$.
Builder proceeds recursively as follows: if there is some $x\in C_1^r\setminus C_2^r$ and $y\in C_2^r\setminus C_1^r$, Builder queries the edge $xy$ and sets $C_{\chi(xy)}^{r+1}=C_{\chi(xy)}^r\cup\{x,y\}$ and $C_{3-\chi(xy)}^{r+1}=C_{3-\chi(xy)}^r$.

Notice that, for every $r$, $C_i^r$ is a connected component in color $i$, so if ever $|C_i^r|=n$ for either $i=1$ or $i=2$, then Builder has located a spanning tree in color $i$.
Set $a(r)=|C_1^r|+|C_2^r|$, so $a(0)=n+1$ and if $a(r)\geq 2n-1$, then it must be the case that $|C_i^r|=n$ for either $i=1$ or $i=2$.
If Builder has not located a monochromatic spanning tree by step $r$, then since $C_1^r\cup C_2^r=V(K_n)$, we must have $a(r+1)=a(r)+1$.
We conclude that Builder can locate a monochromatic spanning tree using at most $(n-1)+(n-2)=2n-3$ queries.
\qed

\paragraph{Painter's strategy.}
Painter will color the first $n-2$ queried edges with color $1$, color the next $n-2$ queried edges with color $2$, and then color the remaining edges arbitrarily.
Since a spanning tree has $n-1$ edges, Builder cannot have located a monochromatic spanning tree within the first $2(n-2)$ queries, and thus must query at least $2n-3$ edges.
\qed
\subsection{Three colors}\label{sec:tree}

\paragraph{Builder's strategy.}
We begin with the key lemma which motivates Builder's strategy.
Recall that $k(n):= {n\over 2}+1$ if $n\equiv 2\pmod 4$ and $k(n):=\bigl\lceil{n\over 2}\bigr\rceil$ otherwise.

\begin{lemma}\label{lem:6cover}
	For $n\geq 3$, suppose that there are $U_1,\dots,U_6\subseteq V(K_n)$ (some of which may be empty) with $E(K_n)=\bigcup_{i=1}^6{U_i\choose 2}$.
	If $U_1\cap U_2=\varnothing$, then $|U_i|\geq k(n)$ for some $i\in[6]$.
\end{lemma}

\begin{proof}
	Assume for the sake of contradiction that $|U_i|<k:=k(n)$ for each $i\in[6]$.
Since $k-1\leq\lfloor n/2\rfloor$, and $U_1\cap U_2=\varnothing$, we can find $A\supseteq U_1$ and $B\supseteq U_2$ with $A\cap B=\varnothing$ and $|A|=\lceil n/2\rceil$ and $|B|=\lfloor n/2\rfloor$.

Consider the edges between $A$ and $B$, denoted by $E[A,B]$.
Since $|U_i|\leq k-1$, by convexity we must have
\[
	\biggl| {U_i\choose 2}\cap E[A,B]\biggr|\leq\biggl\lfloor{(k-1)^2\over 4}\biggr\rfloor,
\]
for every $i$.
However, the edges induced by $U_3,\dots,U_6$ must cover all edges between $A$ and $B$, so we must have
\[
	4\biggl\lfloor{(k-1)^2\over 4}\biggr\rfloor\geq \biggl\lfloor{n^2\over 4}\biggr\rfloor,
\]
contradicting the definition of $k=k(n)$ for any value of $n\pmod 4$.
\end{proof}

Hence, Builder will work to find six subsets of $V(K_n)$, each being contained within a connected component of some color class, wherein each pair of vertices are contained in one of these sets.
The following lemma presents the primary tool in Builder's strategy.

\begin{lemma}\label{lem:compextend}
	Let $\chi\colon E(K_n)\to [3]$ and suppose that Builder has queried some edges giving rise to graphs $C_1,C_2,C_3$ where $C_i$ is the graph formed by the exposed edges of color $i$.

	Suppose that $V_1\subseteq V(C_1)$ and $V_2\subseteq V(C_2)$ are subsets of connected components of $C_1$ and $C_2$, respectively.
	There exists a procedure $\compextend(\chi,V_1,V_2)$ which returns a tuple $(X_1,X_2,X_3)$ of subsets of $V(K_n)$ using at most $2|V_1|+2|V_2|$ additional queries.
	If $C_i^*$ denotes the graph formed by the exposed edges of color $i$ after calling $\compextend(\chi,V_1,V_2)$, then $(X_1,X_2,X_3)$ satisfies:
	\begin{enumerate}
		\item\label{item:tree1} $X_1,X_2,X_3\subseteq V_1\cup V_2$,
		\item\label{item:tree2} $X_1\supseteq V_1$ and $X_2\supseteq V_2$,
		\item\label{item:tree3} for all $i\in[3]$, $X_i$ is a subset of some connected component of $C_i^*$, and
		\item\label{item:tree4} one of the following:
			\begin{enumerate}
				\item $X_1=V_1\cup V_2$,
				\item $X_2=V_1\cup V_2$, or
				\item $X_3\supseteq (V_1\setminus X_2)\cup(V_2\setminus X_1)$.
			\end{enumerate}
	\end{enumerate}
\end{lemma}
\begin{proof} We first define $\compextend$ and then prove the claimed properties.
\begin{algorithmic}[1]
\Procedure{CompExtend}{$\chi,V_1,V_2$}
	\State $X_1\gets V_1$
	\State $X_2\gets V_2$
	\State $X_3\gets \varnothing$
	\Loop
		\If{there is $u\in X_1\setminus(X_2\cup X_3)$ \textbf{and} $v\in(X_2\cap X_3)\setminus X_1$}
			\State \query the edge $uv$
			\State $X_{\chi(uv)}\gets X_{\chi(uv)}\cup\{u,v\}$\label{algl:ex1}
		\ElsIf{there is $u\in X_2\setminus(X_1\cup X_3)$ \textbf{and} $v\in(X_1\cap X_3)\setminus X_2$}
			\State \query the edge $uv$
			\State $X_{\chi(uv)}\gets X_{\chi(uv)}\cup\{u,v\}$\label{algl:ex2}
		\ElsIf{$X_3\subseteq X_1\cap X_2$ \textbf{and} there is $u\in X_1\setminus X_2$ \textbf{and} there is $v\in X_2\setminus X_1$}
			\State \query the edge $uv$
			\If{$\chi(uv)\in\{1,2\}$}
				\State $X_{\chi(uv)}\gets X_{\chi(uv)}\cup\{u,v\}$\label{algl:ex3}
			\Else
				\State $X_3\gets\{u,v\}$\label{algl:restart}
			\EndIf
		\Else\label{algl:end}
			\State\Return $(X_1,X_2,X_3)$
		\EndIf
	\EndLoop
\EndProcedure
\end{algorithmic}
	Items~\eqref{item:tree1} and~\eqref{item:tree2} are straightforward to check.
	To verify Item~\eqref{item:tree3}, notice that in Lines~\ref{algl:ex1},~\ref{algl:ex2} and~\ref{algl:ex3}, if we extended $X_i$ to include some vertex $u$, then this is because there is some $v\in X_i$ with $\chi(uv)=i$.
	In particular, since $V_1,V_2$ were subsets of connected components of $C_1,C_2$, respectively, we have that at every stage, $X_i$ is a subset of a connected component in the currently exposed edges of color $i$.

	Item~\eqref{item:tree4} follows from the observation that this is the only way that the procedure will ever break the loop and return $(X_1,X_2,X_3)$.
	Thus, to verify Item~\eqref{item:tree4}, we simply must show that the procedure does indeed terminate.
	We do so by showing that the procedure terminates after at most $2|V_1|+2|V_2|+1$ loops, implying also that it queried at most $2|V_1|+2|V_2|$ additional edges.

	For an integer $r$, let $X_i^r$ denote the value of $X_i$ in $\compextend(\chi,V_1,V_2)$ \emph{before} running through the loop for the $(r+1)$st time, so $X_1^0=V_1$, $X_2^0=V_2$ and $X_3^0=\varnothing$.
	Notice that we always have $X_1^r\subseteq X_1^{r+1}$ and $X_2^r\subseteq X_2^{r+1}$ and $X_1^r\cup X_2^r=V_1\cup V_2\supseteq X_3^r$.

	Set $a(r):=2|X_1^r|+2|X_2^r|+|X_3^r\setminus(X_1^r\cap X_2^r)|$.
	Certainly $a(0)=2|V_1|+2|V_2|$ and for any $r$,
	\begin{align*}
		a(r) &= 2|X_1^r\cup X_2^r|+2|X_1^r\cap X_2^r|+|X_3^r\setminus(X_1^r\cap X_2^r)|\\
		     &= 2|X_1^r\cup X_2^r|+|X_1^r\cap X_2^r|+|X_3^r\cup(X_1^r\cap X_2^r)|\\
		     &\leq 3|X_1^r\cup X_2^r|+|X_1^r\cap X_2^r|\leq 4|V_1\cup V_2|\leq 4|V_1|+4|V_2|.
	\end{align*}
	We claim that if $\compextend(\chi,V_1,V_2)$ has \emph{not} output $(X_1,X_2,X_3)$ after the $(r+1)$st loop (that is, the procedure will run through the loop for an $(r+2)$nd time), then $a(r+1)\geq a(r)+1$, thus implying the claim.
	We break into the following cases depending on how loop $r+1$ terminates:
	\begin{enumerate}
		\item If the procedure reaches Line~\ref{algl:ex1}, then one of the following occurs:
			\begin{enumerate}
				\item $X_3^{r+1}=X_3^r\cup\{u\}$. Here, since $u\notin X_2^r\cup X_3^r$, we have $|X_3^{r+1}\setminus (X_1^{r+1}\cap X_2^{r+1})|=|X_3^r\setminus(X_1^r\cap X_2^r)|+1$, so $a(r+1)=a(r)+1$.
				\item $|X_i^{r+1}|\geq |X_i^r|+1$ for some $i\in\{1,2\}$. Thus, even though $|X_3^{r+1}\setminus(X_1^{r+1}\cap X_2^{r+1})|$ may decrease by $1$, we still have $a(r+1)\geq a(r)+2-1=a(r)+1$.
			\end{enumerate}
		\item If the procedure reaches Line~\ref{algl:ex2}, then $a(r+1)\geq a(r)+1$ by a symmetric argument.
		\item If the procedure reaches Line~\ref{algl:ex3}, then we have $|X_i^{r+1}|\geq|X_i^r|+1$ for some $i\in\{1,2\}$ and also $|X_3^{r+1}\setminus(X_1^{r+1}\cap X_2^{r+1})|=|X_3^r\setminus(X_1^r\cap X_2^r)|=0$, so $a(r+1)=a(r)+2$.
		\item Finally, if the procedure reaches Line~\ref{algl:restart}, then $X_i^{r+1}=X_i^r$ for $i\in\{1,2\}$ and $X_3^{r+1}=\{u,v\}$. Since $u,v\notin X_1^r\cap X_2^r$, we have $|X_3^{r+1}\setminus(X_1^{r+1}\cap X_2^{r+1})|=2$ while $|X_3^r\setminus(X_1^r\cap X_2^r)|=0$, so again $a(r+1)=a(r)+2$.\qedhere
	\end{enumerate}
\end{proof}

We now have all of the necessary tools describe Builder's full strategy and prove Theorem~\ref{thm:tree}.

\begin{proof}[Proof of Theorem~\ref{thm:tree}]
	Let $\chi\colon E(K_n)\to\{r,g,b\}$ be a 3-coloring.
	Builder begins by choosing $v\in V(K_n)$ arbitrarily and queries all $n-1$ edges incident to $v$.
	Let $R=\{u\in V(K_n)\setminus\{v\}:\chi(vu)=r\}$ and define $G$ and $B$ analogously for colors $g$ and $b$.
	Note that $R,G,B$ are subsets of the vertices of some connected components in the currently exposed edges in colors $r,g,b$, respectively.
	Now, Builder uses $\compextend(\chi,R,G)$, $\compextend(\chi,G,B)$ and $\compextend(\chi,B,R)$ (with the appropriate relabeling of the colors) to find tuples $(R_1,G_2,B_3)$, $(G_1,B_2,R_3)$ and $(B_1,R_2,G_3)$, respectively, as in Lemma~\ref{lem:compextend}.
	This requires at most
	\[
		\bigl(2|R|+2|G|\bigr)+\bigl(2|R|+2|B|\bigr)+\bigl(2|B|+2|G|\bigr)=4\bigl(|R|+|B|+|G|\bigr)=4(n-1)
	\]
	additional queries, thus bringing the total number of queries to at most $5(n-1)$.

	We claim that Builder has located a monochromatic tree on at least $k(n)$ vertices.
	Suppose that $C_r,C_g,C_b$ are the graphs formed by the exposed edges in colors $r,g,b$, respectively.
	If $R_1\cup R_2\cup R_3$ is a subset of a connected component of $C_r$, set $R_1^*=R_1\cup R_2\cup R_3\cup\{v\}$ and $R_2^*=\varnothing$, and otherwise set $R_1^*=R_1\cup R_2\cup\{v\}$ and $R_2^*=R_3$.
	Define $G_i^*$ and $B_i^*$ analogously for $i\in\{1,2\}$.

	In any case, by Lemma~\ref{lem:compextend}, we know that for each $i\in\{1,2\}$, $R_i^*$ is a subset of the vertices of a connected component of $C_r$, $G_i^*$ is a subset of the vertices of a connected component of $C_g$, and $B_i^*$ is a subset of a connected component of $C_b$, thus we need only show that at least one of these sets has size at least $k(n)$.
	We do this by appealing to Lemma~\ref{lem:6cover}.

	By definition, $R_1^*\cap R_2^*=G_1^*\cap G_2^*=B_1^*\cap B_2^*=\varnothing$, so we need only show that
	\[
		E(K_n)={R_1^*\choose 2}\cup{R_2^*\choose 2}\cup{G_1^*\choose 2}\cup{G_2^*\choose 2}\cup{B_1^*\choose 2}\cup{B_2^*\choose 2},
	\]
	i.e.\ every pair of vertices of $K_n$ are contained together in one of these six sets.
	Let $x,y\in V(K_n)$ be two distinct vertices.
	If, say, $x=v$, then $y\in R\cup B\cup G$, and so $x,y$ are contained together in one of $R_1^*,G_1^*,B_1^*$.
	If we have $x,y\in R$, then certainly $x,y\in R_1^*$, and similarly if $x,y\in G$ or $x,y\in B$.

	Thus, suppose that, without loss of generality, $x\in R,y\in G$.
	Certainly $x\in R_1\subseteq R_1^*$ and $y\in G_2\subseteq G_1^*$, so suppose that $x\notin G_1^*$ and $y\notin R_1^*$.
	Then by Item~\ref{item:tree4} in Lemma~\ref{lem:compextend}, we must have $B_3\supseteq (R\setminus G_2)\cup(G\setminus R_1)$, so $x,y\in B_3$.
	Hence, either $x,y\in B_1^*$ or $x,y\in B_2^*$.
\end{proof}

\paragraph{Painter's strategy.}
Painter's strategy here is very similar to her strategy in the proof of Theorem~\ref{thm:matchingtight}.

\begin{proof}[Proof of Theorem~\ref{thm:treetight}]
    Since $\oram_3(\cal T_{k(n)+1};n)\geq\oram'_3(\cal T_{k(n)+1};n)$, it suffices to give only a strategy for the cornering game.

	Let $\chi$ be a $3$-coloring of $E(K_n)$ formed by starting with a proper edge-coloring of $K_4$ using three colors, ``blowing up'' each vertex into a cluster of size roughly $n/4$, and then coloring the edges within the clusters arbitrarily.
	Formally speaking, start with a partition $[n]=V_1\sqcup V_2\sqcup V_3\sqcup V_4$ with $|V_i|\in\bigl\{\lfloor n/4\rfloor,\lceil n/4\rceil\bigr\}$ for all $i$, and let $\chi$ be the $3$-coloring of $E(K_n)$ given by
	\begin{itemize}
		\item $\chi(e)=r$ if $e\in E[V_1,V_2]\cup E[V_3,V_4]$,
		\item $\chi(e)=b$ if $e\in E[V_1,V_3]\cup E[V_2,V_4]$,
		\item $\chi(e)=g$ if $e\in E[V_1,V_4]\cup E[V_2,V_3]$, and
		\item $\chi(e)$ is arbitrary otherwise.
	\end{itemize}

	It is straightforward to verify that $\chi$ does not contain a monochromatic tree on $k(n)+1$ vertices.

	For an edge $e$ and color $c\in\{r,b,g\}$, let $\chi_{e,c}$ denote the coloring where $e$ gets color $c$ and every other edge is colored as in $\chi$.
	Notice that if $e$ is any edge not completely contained in some $V_i$, then $\chi_{e,c}$ actually contains a spanning tree in color $c$ whenever $c\neq \chi(e)$.

	When Builder queries an edge $e$, Painter colors $e$ as in $\chi$, unless $e$ is the last unexposed edge which is not completely contained in one of the $V_i$'s.
	In this situation, Painter gives $e$ either color $c_1$ or $c_2$ where $c_1,c_2\neq\chi(e)$.

	Thus, Painter's coloring will always contain a monochromatic tree on $k(n)+1$ vertices (in fact, it will always contain a monochromatic spanning tree), but Builder must query every edge not completely contained in some $V_i$ to determine which color contains said tree.
	As such, Builder must query at least $\sum_{i\neq j\in[4]}|V_i||V_j|\geq 6\bigl\lfloor{n\over 4}\bigr\rfloor^2$	edges to do so.
\end{proof}

\section{Remarks}\label{sec:remarks}
We begin by showing that $\oram'$ is indeed a generalization of $\oram$.
\begin{lemma}\label{lem:locvcor}
    Fix $t\geq 2$ and let $\cal G_1,\dots,\cal G_t$ be families of graphs.
    If $n\geq\ram(\cal G_1,\dots,\cal G_t)$, then $\oram'(\cal G_1,\dots,\cal G_t;n)=\oram(\cal G_1,\dots,\cal G_t;n)$.
\end{lemma}
\begin{proof}
    We need only show the inequality $\oram(\cal G_1,\dots,\cal G_t;n)\leq\oram'(\cal G_1,\dots,\cal G_t;n)$.
    By the assumption on $n$, any strategy for Painter in the locating game is also a strategy for Painter in the cornering game.

    Suppose that Builder has exposed some collection of edges, giving rise to color classes $C_1,\dots,C_t$, and that Builder has won the cornering game by determining, without loss of generality, that color $1$ must contain a copy of some $G_1\in\cal G_1$.
    Consider the coloring $\chi\colon E(K_n)\to[t]$ defined by $\chi(e)=i$ if $e\in C_i$ and $\chi(e)=t$ otherwise.
    Since $\chi$ could be Painter's hidden coloring and $C_1=\chi^{-1}(1)$, it must be the case that $C_1$ contains some $G_1\in\cal G_1$.
    Therefore, Builder has additionally won the locating game.
\end{proof}
\medskip

We showed in Theorem~\ref{thm:matching} that $\oram_t(rK_2;n)\leq O(n\log t)$ whenever $n\geq(t+1)r-t+1$.
While we have already remarked that the precise bound in Theorem~\ref{thm:matching} is not necessarily tight (Remark~\ref{rm:nottight}), we wonder if it is tight up to a constant:
\begin{question}\label{q:matching}
	Fix $t,r\geq 2$ and set $n=(t+1)r-t+1$. Is $\oram_t(rK_2;n)\geq\Omega(n\log t)$?
\end{question}
Notice that, by Lemma~\ref{lem:locvcor}, it does not matter whether we consider $\oram$ or $\oram'$ in this question.

We additionally do not know the correct answer even in the case of only $2$ colors.
In this situation, since $n=3r-1\not\equiv 0\pmod 3$, we have $\oram_2(rK_2;n)\leq n-1$ (Remark~\ref{rm:constimprove}).
\begin{question}
	For $r\geq 1$ and $n=3r-1$, is it the case that $\oram_2(rK_2;n)=n-1$?
\end{question}

\medskip

Turning now to trees, we showed in Theorem~\ref{thm:tree} that $\oram_3(\cal T_{k(n)};n)\leq 5(n-1)$.
We wonder if this result can be improved:
\begin{question}
	What is the smallest $c$ such that in any $3$-coloring of $E(K_n)$, Builder can locate a monochromatic tree on $k(n)$ vertices by querying at most $\bigl(c+o(1)\bigr)n$ edges?
\end{question}
Again, thanks to Lemma~\ref{lem:locvcor}, it does not matter whether we consider $\oram$ or $\oram'$ here.
Theorem~\ref{thm:tree} shows that $c\leq 5$; furthermore, certainly $c\geq 3/2$ since Painter can simply color the first $k(n)-1$ edges red, the next $k(n)-1$ edges blue, and the next $k(n)-1$ edges green.

Theorems~\ref{thm:2ramtree} and~\ref{thm:ramtree} were extended by Gy\'arf\'as~\cite{G77} to show that any $t$-coloring of $E(K_n)$ must contain a monochromatic tree on at least ${n\over t-1}$ vertices. F\"uredi~\cite{furedi1989covering} showed that this bound can be improved slightly in the case where an affine plane of order $t-1$ does not exist.
These suggest the natural question:
\begin{question}
	What is the least number of queries necessary for Builder to locate a monochromatic tree on at least ${n\over t-1}$ vertices in a $t$-coloring of $E(K_n)$?
\end{question}
We suspect that $\oram_t(\cal T_{n/(t-1)};n)\leq O_t(n)$, which was the case for both $2$ and $3$ colors.

\paragraph{Acknowledgments.} We are grateful to the anonymous referee, whose comments helped significantly improve the presentation of our results.

\bibliographystyle{abbrv}
\bibliography{references}

\end{document}